\documentclass{svmult}

\usepackage{makeidx}         
\usepackage{graphicx}        
\usepackage{multicol}        
\usepackage[bottom]{footmisc}
\usepackage{url}

\usepackage{epsfig}

\renewcommand{\PackageWarningNoLine}[2]{}
\usepackage{amsmath}
\usepackage{amsfonts}

\begin{document}

\author{
Juan Galvis\inst{1}
\and
Eric Chung\inst{3}
\and
Yalchin Efendiev\inst{2}
\and
Wing Tat Leung\inst{2}
}
\authorrunning{}

\institute{%
$^1$Departamento de Matem\'aticas, Universidad Nacional de Colombia, Bogot\'a, Colombia.
$^2$Department of Mathematics, Texas A\&M University,
College Station, TX 77843-3368, USA.
$^3$Department of Mathematics, The Chinese University of Hong Kong, Hong Kong SAR.}

\title*{On overlapping domain decomposition methods for high-contrast multiscale problems}

\titlerunning{On overlapping DD for high-contrast multiscale problems}
\authorrunning{Galvis, Chung, Efendiev \& Leung.}


\maketitle
\section{Summary}
\vspace{-0.5cm}
We review some important ideas in the design and analysis of robust overlapping domain decomposition algorithms for high-contrast multiscale problems and propose a domain decomposition method better performance in terms of the number of iterations. The main novelty of our approaches is the construction of coarse spaces, which are computed using spectral information of local bilinear forms. We present several approaches to incorporate the spectral information into the coarse problem in order to obtain minimal coarse space dimension. We show that using these coarse spaces, we can obtain a domain decomposition preconditioner with the condition number independent of contrast and small scales. To minimize further the number of iterations until convergence, we use this minimal dimensional coarse spaces in a construction combining them with large overlap local problems that take advantage of the possibility of localizing global fields orthogonal to the coarse space. We obtain a condition number close to 1 for the new method. We discuss possible drawbacks and further extensions.\vspace{-0.5cm}

\section{High-contrast problems. Introduction}
\vspace{-0.5cm}
The methods and algorithms,
discussed in the paper,
can be applied to various PDEs, even though we
will focus on Darcy flow equations.
Given $ { D}\subset \mathbb{R}^2$, $ {   f:D \to \mathbb{R}}$ ,  and $ {   g:\partial D\to \mathbb{R}}$, find
$u:D\to \mathbb{R}$ such that
\[
{\partial\over\partial x_i} \left( \kappa_{ij}  {\partial u\over\partial x_j} \right) = f
\]
with a suitable boundary condition, for instance $u=0$ on $\partial D$.
The coefficient  $\kappa_{ij}(x)=\kappa(x) \delta_{ij}$
represents the permeability of the porous media $D$.
We focus on two-levels
overlapping domain decomposition and use
local spectral information  in constructing ``minimal'' dimensional coarse spaces (MDCS).
After some review on constructing MDCS and their use in overlapping
domain decomposition preconditioners, we present an approach, which uses
MDCS to minimize the condition number to a condition number closer to 1. This approach
requires a large overlap (when comparted to coarse-grid size) and, thus, is more efficient for
small size coarse grids.
We present the numerical results
and state our main theoretical result.
We assume that there exists $\kappa_{\min}$ and $\kappa_{\max}$ with $
0< \kappa_{\min}\leq \kappa(x)\leq \kappa_{\max}$ for all  $x\in D$.
 {\bf The coefficient $\kappa$ has a   multiscale structure} (significant local variations of $\kappa$ occur across $D$ at
different scales). We also assume that {the coefficient $\kappa$ is a  high-contrast coefficient} (the constrast is $ \eta={\kappa_{\max}}/{\kappa_{\min}}$). We  assume that $\eta$ is large compared to the coarse-grid size.

It is well known that  performance of numerical methods for
 high-contrast multiscale problems depends on $\eta$
 and local variations of $\kappa$ across $D$.
 For classical finite element methods,
 the condition to obtain good approximation results is that the finite
 element mesh has to be fine enough
 to resolve the variations of the coefficient $\kappa$.
 Under these conditions, finite element approximation leads
 to  the solution of very large (sparse)
 ill-conditioned problems (with the condition number scaling with $h^{-2}$ and $\eta$). Therefore, the
performance of solvers  depends on $\eta$
and local variations of $\kappa$ across $D$.
This was observed in  several works, e.g., \cite{ge09_1, Graham1, aarnes}\footnote{Due to the page limitation, only a few references are cited throughout.}.

Let $\mathcal{T}^h$ be a triangulation of the domain $D$, where
$h$ is the size of typical element. We consider only the case
of discretization  by the classical
finite element method $V=P_1(\mathcal{T}^h)$ of piecewise (bi)linear functions. Other discretizations
can also be considered.  The application of the finite element discretization
leads to the solution of a very large ill-conditioned system
 $
 Ax=b,
$
where   $A$ is roughly of size  $h^{-2}$  and the condition number of $A$ scales
with $\eta$ and $h^{-2}$. In general, the main goal is to obtain an  efficient good
approximation of solution $u$.
The two main solution strategies  are:
{\bf 1. Choose $h$ sufficiently small and implement an iterative method}.
It is important to implement a preconditioner
$M^{-1}$ to solve $M^{-1}Au=M^{-1}b$. Then, it is important to have the condition number of
$M^{-1}A$ to be small and bounded independently of physical parameters, e.g.,
$\eta$ and the
multiscale structure of $\kappa$.
{\bf 2. Solve a smaller dimensional linear system}
($\mathcal{T}^H$ with {$H>h$}\footnote{The coarse mesh does not necessarily resolve
  all the variations of $\kappa$} ) so that computations of solutions   can be done efficiently.
This usually involves the
construction of  a downscaling operator $R_0$ (from the coarse-scale  to fine-scale $v_0\mapsto v$)
and an upscaling operator (from fine-scale to coarse-scale, $v\mapsto v_0$) (or similar operators).
Using these operators, the linear system $Au=b$ becomes  a coarse linear system
$A_0 u_0=b_0$ so that  $R_0u_0$ or functionals of it can be computed. The main
goal of this approach it to obtain a sub-grid capturing such that
 $|| u- R_0u_0||$  is small.


The rest of the paper will focus on the design of
overlapping domain decomposition methods by constructing appropriate coarse spaces.
First, we will review existing results, which construct minimal dimensional coarse spaces,
such that the condition number of resulting preconditioner is independent of $\eta$.
These coarse spaces use local spectral problems to extract the information, which can not be localized.
This information is related to high-conductivity channels, which connect coarse-grid boundaries
and important in domain decomposition preconditioners and multiscale simulations.
Next, using these MDCS and oversampling ideas, we present a ``hybrid'' domain decomposition approach
with a condition number close to 1 by appropriately selecting the oversampling size (i.e., overlapping size). We state
our main result, discuss some limitations, and show a numerical example.
We compare the results to some existing contrast-independent
preconditioners.


\vspace{-0.5cm}
\section{Classical overlapping methods. Brief review}
\vspace{-0.5cm}
We start with a non-overlapping decomposition $\{ D_i\}_{i=1}^{N_S}$ of the domain $D$
and obtain an overlapping decomposition $\{ D_i'\}_{i=1}^{N_S}$
by adding a layer of  width $\delta$ around each non-overlapping subdomain.
Let $A_j$ be the  Dirichlet matrix corresponding to the overlapping subdomain $D_j'$.
The one level method solves
$M_1^{-1}A=M_1^{-1}b$ with
$M^{-1}_1\sum_{j=1}^{N_S} R_j({A}_j)^{-1} R_j^T$ and the operators
$R_j^T$, $j=1,\dots,N_S$, being the restriction to overlapping subdomain $D_j'$ operator and with
the  $R_j$ being the extension by zero (outside $D_j'$) operator.
We have the bound $\mbox{Cond}(M^{-1}_1A)\leq C\left(1+{1}/{\delta H}\right)$.
For high-contrast multiscale problems, it is known
that $\displaystyle C \asymp \eta$.

Next, we introduce a coarse space, that is, a subspace $V_0\subset V$ of small dimension
(when compared to the fine-grid finite element space V).
We consider $A_0$ as the matrix form of the discretization of the equation related to subspace
$V_0$. For simplicity of the presentation, let  $A_0$  be the Galerkin projection of
$A$ on the subspace $V_0$. That is
$A_0=R_0AR_0^T$,
where $R_0$ is a downscaling operator
that converts coarse-space coordinates into
fine-grid space coordinates. The two-levels preconditioner
uses the coarse space and it is defined by $M^{-1}_2=  R_0A_0^{-1}R_0^T+
\sum_{j=1}^{N_S} R_j ({A}_j)^{-1} R_j^T= R_0A_0^{-1}R_0^T+M_1^{-1}$.
It is known that $\mbox{Cond}(M^{-1}A)\preceq
\eta
\left(1+{H}/{\delta}\right).$ The classical two-levels method is robust with respect to the number of subdomains but it is not robust with respect to $\eta$.
The condition number  estimates use Poincar\'e
inequality and a small overlap trick; \cite{tw}. Without  small overlap trick
$\mbox{Cond}(M^{-1}A)\preceq \eta(1+H^2/\delta^2)$.

There were several works addressing the performance of classical domain
decomposition algorithms for high-contrast problems. Many of these works
considered simplified multiscale structures\footnote{These works usually assume
  some alignment between the coefficient heterogeneities and
  the initial non-overlapping decomposition}, see e.g., \cite{tw}
for some works by O. Widlund and his
collaborators.
We also mention  the
works by Sarkis  and his collaborators, where they introduce the assumption of
quasi-monotonicity \cite{MR1367653}.
Sarkis also introduced the idea of using ``extra'' or additional basis
functions as well as techniques that construct the coarse spaces using the
overlapping decomposition (and not related to a coarse mesh); \cite{MR2099424}. Scheichl and  Graham \cite{Graham1} and  Hou and Aarnes \cite{aarnes}, started a
systematic study of the performance of classical overlaping domain decomposition methods for
high-contrast problems. In their works, they used coarse spaces constructed using a
coarse grid and special basis functions from the family of  multiscale finite element
methods.  These  authors designed two-levels
domain decomposition methods that were robust
(with respect to $\eta$) for special
multiscale structures. None of the results available
in the literature
(before the method in papers \cite{ge09_1,ge09_1reduceddim} was introduced)
were robust for
a coefficient not-aligned with the
construction of the coarse space (i.e., not aligned either with the non-overlapping
decomposion or the coarse mesh if any), i.e.,
the condition number of the
resulting preconditioner
is independent of $\eta$ for general multiscale coefficients.

 \vspace{-0.5cm}

\section{Stable decomposition and eigenvalue problem. Review}
\vspace{-0.5cm}
A main tool in obtaining condition number bounds is the construction of a stable decomposition of a global
field. That is, if  for all $v\in V=P^1_0(D,\mathcal{T}^h)$
there exists a decomposition $
 v=v_0+\sum_{j=1}^{N_S} v_j$
 with $v_0\in V_0$ and $ v_j\in V_j=P^1_0(D_j', \mathcal{T}^h)$,
 $j=1,\dots, N$,
and
\[
\int_{D} \kappa |\nabla v_0|^2+\sum_{j=1}^{N_S}
\int_{D_j'}\kappa |\nabla v_j|^2 \leq C_0^2\int_{D}\kappa |\nabla v|^2
\]
for $C_0>0$. Then,
$\mbox{cond}(M_2^{-1}A)\leq c(\mathcal{T}^h,\mathcal{T}^H) C_0^2$.
  Existence of a suitable coarse interpolation
$I_0:V\to V_0=\mbox{span}\{\Phi\}$   implies the stable
  decomposition above. Usually such
  stable decomposition is constructed as follows.

For the coarse part of the stable decomposition, we
introduce a partition of unity  $\{\chi_i\}$ subordinated to the coarse mesh
(supp $\chi_i\subset \omega_i$ where $\omega_i$ is
the coarse-block neighborhood of the coarse-node $x_i$).
We  begin by restricting the global field $v$ to  $\omega_i$.
For each coarse node neighborhood $\omega_i$, we
identify local field that will contribute to the
coarse space
{$  I_0^{\omega_i}v $} so that the
coarse space will be defined as  $V_0= \mbox{Span}\{ \chi_i  I_0^{\omega_i}v\}$.
We assemble a coarse field as
$v_0=I_0v=\sum_{i=1}^{N_S} \chi_i (I_0^{\omega_i}v )$. Note that
in each block $  v-v_0=\sum_{i\in K} \chi_i (v-  I_0^{\omega_i}v )$.

For the local parts of the stable decomposition, we
introduce a partition of unity  $\{\xi_j\}$ subordinated to the non-overlapping
decomposition (supp $\xi_j\subset  D_j'$).  The local part of the stable decomposition is defined by
 $  v_j=\xi_j (v-v_0)$. For instance, to bound the energy of  $v_j$, we have in each coarse-block $K$,
\begin{eqnarray*}
&&\int_{K} \kappa|\nabla v_j|^2
 \preceq \int_{K} \kappa|\nabla  \xi_j\left(\sum_{i\in K} \chi_i (v-  I_0^{\omega_i}v )\right)|^2\\
&\preceq&\sum_{i\in K}\int_{K} \kappa(\xi_j\chi_i)^2|\nabla (v-I_0^{\omega_i}v)|^2
+{   \sum_{x_i\in K} \int_{K} \kappa|\nabla (\xi_j\chi_i)|^2| v-
I_0^{\omega_i}v  |^2}.
\end{eqnarray*}

Adding up over $K$, we obtain,
\begin{eqnarray*}
\int_{D_j' } \kappa|\nabla v_j|^2 \preceq\sum_{i\in D_j'}\int_{D_j' } \kappa(\xi_j\chi_i)^2|\nabla (v-I_0^{\omega_i}v)|^2\\\
+{   \sum_{x_i\in \omega_j} \int_{D_j'} \kappa|\nabla (\xi_j\chi_i)|^2| v-
I_0^{\omega_i}v  |^2}\\
\end{eqnarray*}
and we would like to bound the last term by $      C{   \int_{D_j' } \kappa |\nabla v|^2}$.

 For simplicity of our presentation, we consider the case when
 the coarse elements coincide with the non-overlapping decomposition subdomains.
 That is, $D_j'=\omega_j$. In this case, we can replace $\xi$ by $\chi$ and
 replace $\nabla (\chi^2)$ by $\nabla \chi$ so that we need to bound
 $\sum_{x_i\in \omega_j} \int_{\omega_j} \kappa|\nabla   \chi_i |^2| v-
I_0^{\omega_i}v  |^2$. We refer to this design as {\bf coarse-grid based}.

\begin{remark}[General case and overlapping decomposition based design]
  Similar analysis holds in the case when there is no coarse-grid and the coarse space
  is spanned by the partition of unity  $\{\xi_j\}$. We can replace $\chi$ by $\xi$ and
  $\nabla (\xi^2)$ by $\nabla \xi$.
 In general these two partitions are not  related (see Sec. \ref{secabs}).
\end{remark}

We now review the three main arguments to complete the  required bound: 1) Poincar\'e inequality. 2) $L^\infty$ estimates. 3)    Eigenvalue problem.

{\bf 1. Poincar\'e inequality:} Classical analysis  uses Poincar\'e inequality to obtain the required bound above.
That is, the inequality
 $ \frac{1}{H^2} \int_{\omega} (v-\bar{v})^2 \leq C \int_\omega |\nabla v|^2$ to obtain
${   \sum_{x_i\in \omega_j} \int_{\omega_j} \kappa|\nabla   \chi_i |^2| v-
I_0^{\omega_i}v  |^2}\preceq
{  \frac{1}{H^2} \int_{\omega_i}
\kappa| v- I_0^{\omega_i}v  |^2} \preceq
      C{   \int_{\omega_i} \kappa |\nabla v|^2}.$
      In this case,  $I_0^{\omega_i}v$ is the average of $v$ on the subdomain.
      For the case of high-contrast coefficients,  $C$ depends on
      $\eta$, in general. For quasi-monotonic like coefficient  it can be obtained that $C$ is independent of the contrast \cite{MR1367653}.
We also mention  \cite{ge09_1} for the case { \bf locally connected high-contrast  region}. In this
case  $I_0^{\omega_i}v$ is a weighted average. From the argument given in \cite{ge09_1}, it
was clear that when the high-contrast regions break across the domain, defining only one average was not
enough to obtain contrast independent constant in the Poincar\'e inequality.

\begin{figure}
\centering
\psfig{figure=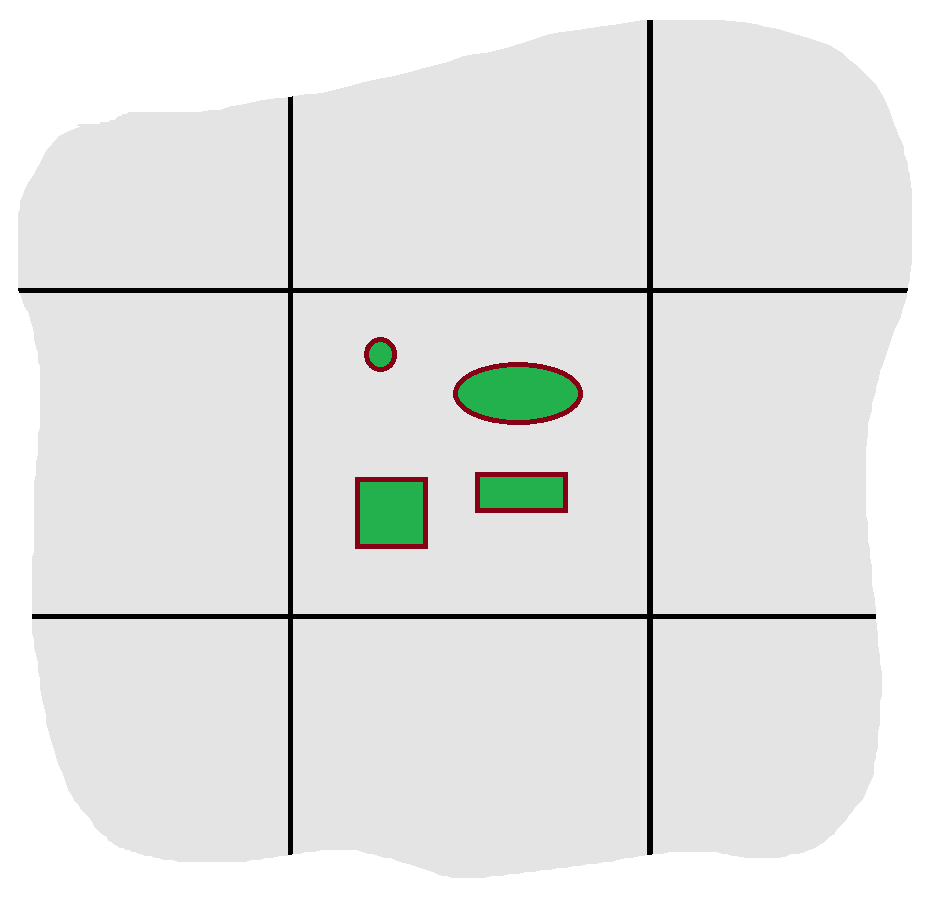,height=3cm,width=3cm,angle=0}
\psfig{figure=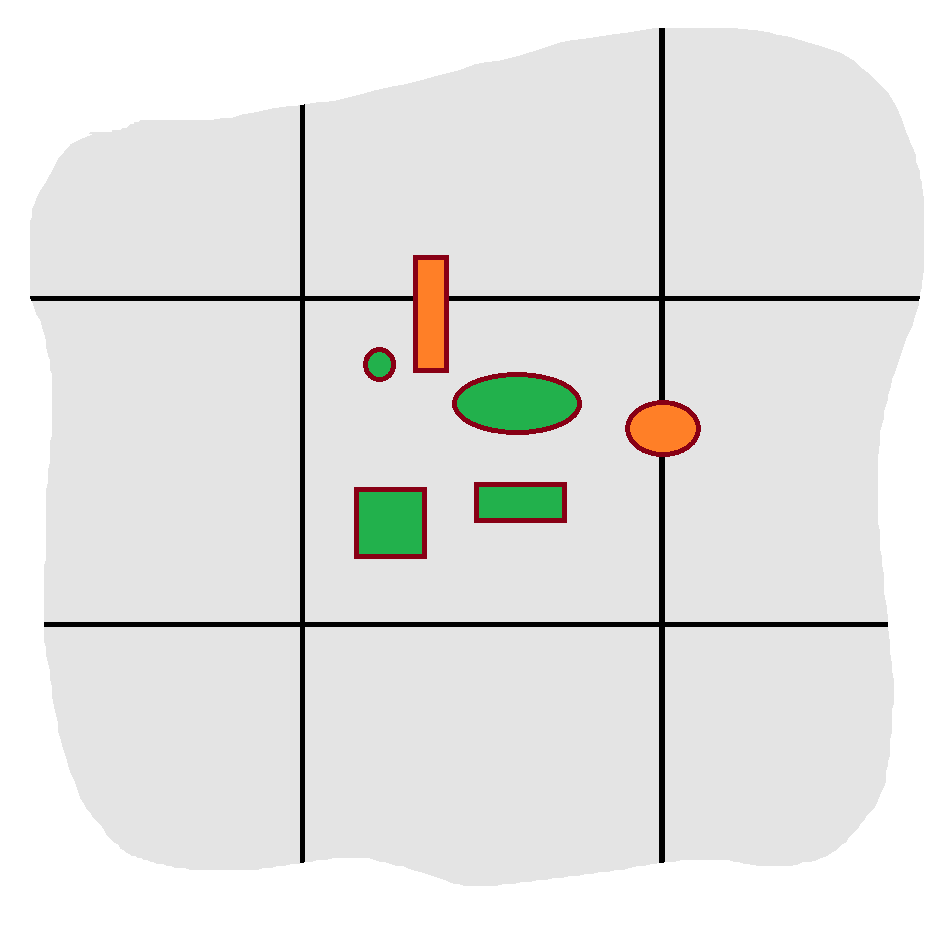,height=3cm,width=3cm,angle=0}
\psfig{figure=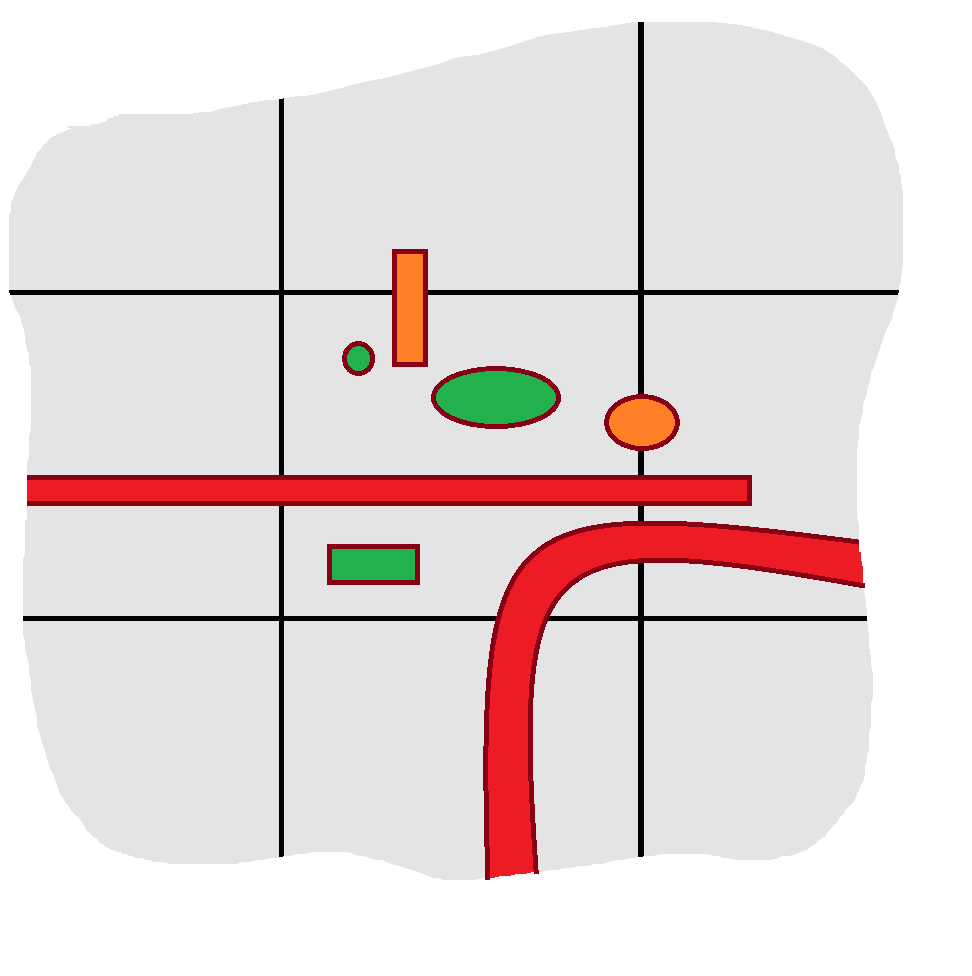,height=3cm,width=3cm,angle=0}
\caption{Examples o multiscale coefficients with interior high-contrast inclusions (left), boundary inclusions (center) and long channels(right).}\label{inclusions}
\end{figure}

{\bf 2. $L^\infty$ estimates:} Other idea is to use an $L^\infty$ estimate of the form
\begin{eqnarray*}
{   \sum_{x_i\in K}\int_{\omega_i} {  \kappa|\nabla \chi_i|^2}| v-
I_0^{\omega_i}v  |^2}
\preceq
   \sum_{x_i\in K}{  || \kappa|\nabla \chi_i|^2 ||_\infty} {   \int_{\omega_i} | v- I_0^{\omega_i}v  |^2}.
\end{eqnarray*}
The idea in \cite{Graham1,aarnes}
was then to construct partition  of unity such that $|| \kappa|\nabla \chi_i|^2 ||_\infty$ is bounded independent of the contrast and then to use classical
 Poincar\'e inequality estimates.
 Instead of minimizing the $L^\infty$,  one can intuitively  try to minimize $\int_K \kappa|\nabla \chi_i|^2 $.
This works well when the multicale structure of the coefficient is confined withing the coarse blocks. For instance, for a coefficient and coarse-grid as depicted in Figure \ref{inclusions} (left picture), we have that a two-level
domain decomposition method can be proven to be robust with respect the value of the coefficient inside
the inclusions. In fact, the coarse space spanned by classical multiscale basis functions with linear boundary
conditions ($-\mbox{div}(\kappa \nabla \chi_i)=0$ in $K$ and linear on each edge of $\partial K$) is sufficient
and the above proof works.
Now consider the coefficient in Figure \ref{inclusions} (center picture).
For such cases, the boundary condition of the basis
functions is important. In these cases,
basis functions can be constructed such that
the above argument can be carried on. Here,
we can use multiscale basis functions with oscillatory boundary
condition in its construction\footnote{We can include  constructions of boundary conditions using $1D$ solution of the problem along the edges. Other choices include basis functions constructed using oversampling regions, energy minimizing partition of unity (global), constructions using limited global information (global), etc.}.

\vspace{-0.5cm}

For the coefficient in Figure \ref{inclusions}, right figure, the argument above using $L^\infty$  cannot be carried out
unless we can work with larger support basis functions (as large as to include the high-contrast channels  of
the coefficient). If the support of the coarse basis function does not include the high-contrast region, then
 $|| \kappa|\nabla \chi_i|^2 ||_\infty$ increases with the contrast leading to non-robust two-level domain
 decomposition methods.

{\bf 3. Eigenvalue problem.} We can write
$\displaystyle{   \sum_{x_j \in \omega_i}\int_{\omega_i} \kappa   |\nabla \chi_j|^2  |v-
I_0^{\omega_i}v |^2}
\preceq
{ \frac{1}{H^2} \int_{\omega_i}
\kappa|(v- I_0^{\omega_i}v )|^2}\preceq C
{  \int_{\omega_i} \kappa |\nabla v|^2}$,
where we need to justify the last inequality with constant independent of the
contrast. The idea is then to consider the Rayleigh quotient,
  $$\mathcal{Q}(v):=\frac{\int_{\omega_i} \kappa |\nabla v|^2}{\int_{\omega_i} \kappa| v  |^2}$$ with
   $v\in P^1(\omega_i)$.
This quotient is related to an eigenvalue problem and we can define
 $ I_0^{\omega_i}v $ to be  the projection on low modes of this quotient  on $\omega_i$.
The associated eigenproblem is given by $-\mbox{div}(k(x)\nabla \psi_\ell)=\lambda_\ell k(x)\psi_i$ in $\omega_i$
with homogeneous Neumann boundary condition for floating subdomains and a mixed homogeneous Neumann-Dirichlet condition for subdomains that touch the boundary.
It turns out that the low part of the spectrum can be written as
$\lambda_1\leq \lambda_2\leq ...\leq \lambda_L$ $<\lambda_{L+1}\leq ...$
where   $\lambda_1 , ...,  \lambda_L$ are small, asymptotically vanishing eigenvalues and $\lambda_L$ can be bounded below independently of the contrast. After identifying the local field $I_0^{\omega_i}v$,
we then define the coarse space as
$V_0=Span\{ I^h\chi_i \psi_j^{\omega_i} \}= Span\{ \Phi_i \}.$

{\bf Eigenvalue problem with multiscale partition of unity.}
Instead of the argument presented earlier,
we can include the gradient of the partition of unity in the bounds (somehow similar to the ideas of $L^\infty$ bounds). We then need the following chain of inequalities,
$\displaystyle {   \int_{\omega_i}\underbrace{
\left( \sum_{x_j \in \omega_i}  \kappa|\nabla \chi_j|^2\right)} _{\displaystyle :=H^{-2}\widetilde{\kappa}}|v-
I_0^{\omega_i}v |^2}
=
{  \frac{1}{H^2} \int_{\omega_i}
\widetilde{\kappa}|v- I_0^{\omega_i}v )|^2}
\preceq
{   \int_{\omega_i} \kappa |\nabla v|^2}.$
Here we have to consider  Rayleigh quotient $\mathcal{Q}_{ms}(v):=\frac{\int_{\omega_i}  \kappa  |\nabla v|^2}{\int_{\omega_i} \widetilde{\kappa}| v  |^2}$,
   $v\in P^1(\omega_i)$  and define $ I_0^{\omega_i}v $ as projection on low modes. Additional modes ``complement'' the initial space spanned by the partition of unity so that the resulting coarse space leads to robust  methods with minimal dimension coarse spaces; \cite{ge09_1reduceddim}.\\

If we
consider the two-level method with the (multiscale) spectral coarse space presented before, then
\begin{equation}\label{b1}
\mbox{cond}(M^{-1}A)\leq C(1+(H/\delta)^2),
\end{equation}
 where $C$ is independent of the contrast
if enough eigenfunctions in each
node neighborhood are selected for the construction
of the coarse spaces. The constant  $C$  and the resulting coarse-space dimension depend on the partition
of unity (initial coarse-grid representation) used.

\vspace{-0.5cm}
\subsection{Abstract problem eigenvalue problems}\label{secabs}
\vspace{-0.5cm}
We consider an abstract variational problem, where the global bilinear form is obtained by
assembling local bilinear forms. That is
$a(u,v)=\sum_K a_K(R_Ku,R_Kv)$, where
$a_K(u,v)$ is a bilinear form acting on functions with supports
being the coarse block $K$. Define the  subdomain bilinear form $a_{\omega_i}(u,v)=\sum_{K\subset \omega_i}a_K(u,v)$.
We consider the abstract problem
$$a(u,v)=F(v) \quad \mbox{ for all } v\in V.$$
We
introduce $\{\chi_j\}$, a partition of unity subordianted to coarse-mesh blocks and
$\{\xi_i\}$ a partition of unity subordianted to overlapping decomposition
(not necessarily related in this subsection). We also define  the
``Mass'' bilinear form (or energy of cut-off) $m_{\omega_i}$ and
the Rayleigh quotient  $\mathcal{Q}_{abs}$ by
 \[\displaystyle m_{\omega_i} (v,v)  :=  \sum_{j\in \omega_i} a(\xi_i\chi_j v,\xi_i\chi_j v)\quad \mbox{ and } \quad \mathcal{Q}_{abs}(v):=\displaystyle \frac{a_{\omega_i} (v,v) }{m_{\omega_i} (v,v) }.
\]

For the Darcy problem,  we have $ m_{\omega_i} (v,v)
= \sum_{j\in \omega_i} \int_{\omega_i} \kappa |\nabla(\xi_i\chi_j v )|^2
\preceq \int_{\omega_i}\widetilde{ \kappa} |v|^2.$
The same analysis can be done by replacing the partition of unity function by  partition of degree of freedom (PDoF). Let  $\{\pmb\chi_j\}$ be PDoF subordianted to coarse mesh neighborhood and
$\{\pmb\xi_i\}$ be PDoF subordianted to overlapping decomposition. As before, we
define the cut-off bilinear form and quotient,
 \[\displaystyle m_{\omega_i} (v,v)  :=  \sum_{j\in \omega_i} a(\pmb\xi_i\pmb\chi_j v,\pmb\xi_i\pmb\chi_j v) \quad \mbox{ and } \quad \mathcal{Q}_{abs2}(v):=\displaystyle \frac{a_{\omega_i} (v,v) }{m_{\omega_i} (v,v) }.\]
 The previous construction alows applying the same design recursively and
 therefore to use the same ideas in a multilevel method. See \cite{eglw11}.

 \subsection{Generalized Multiscale Finite Element Method (GMsFEM) eigenvalue problem}

 We can consider the
 Rayleigh quotients presented before only in a suitable subspace
that allows a good approximation of low modes. We call these subspace the  snapshot spaces. Denote by $W_i$ the
snapshot space corresponding to subdomain  $\omega_i$, then we consider the Rayleigh quotient,
$\displaystyle \mathcal{Q}_{gm}(v):= \frac{a_{\omega_i} (v,v) }{m_{\omega_i} (v,v) } \quad
\mbox{with  $v\in W_i.$}$
The snapshot space can be obtained  by dimension reduction techniques or similar computations.
See \cite{egh12, calo2016randomized}. For example, we can consider the following simple example.
In each subdomain $\omega_i$, $i=1,\dots, N_S$:
(1) Generate forcing terms $f_1,f_2,\dots,f_M$ randomly ($\int _{\omega_i}f_\ell=0$);
(2) Compute the local solutions
$ -\mbox{div}(\kappa \nabla u_\ell )=f_\ell $
with homogeneous Neumann boundary condition; (3) Generate $W_i=\mbox{span}\{ u_\ell \}\cup\{1\}$;
(4) Consider $\mathcal{Q}_{gm}$ with $W_i$ in 3 and compute important modes.\\
In Table 1, we see the results of using the local eigenvalue problem versus using the GMsFEM eigenvalue problem.
\begin{table}[htb]
\centering
\footnotesize
\begin{tabular}{|l|r|r|r|r|r|r|}\hline
$\eta$ & MS&  Full  &  8 rand. &
15 rand \\\hline
  $10^6$    & 209   &    35  & 37  & 37\\
  $10^9$    & 346   &    38 & 44  &38\\
\hline\hline\hline
\end{tabular}
\caption{PCG iterations for different values  $\eta$.
Here $H=1/10$ with $h=1/200$.
We use the GMsFEM eigevalue problem with $W_i=V_i$ (full local fine-grid space),  column 2; $W_i$ spanned by $8$ random samples, column 4, and $W_i$ spanned by 15 samples, column 5. }
\label{tab:perm1}
\end{table}

\vspace{-0.5cm}

\section{Constrained coarse spaces, large overlaps, and DD}
\vspace{-0.5cm}
In this section,
we introduce a hybrid overlapping domain decomposition preconditioner. We use the coarse spaces constructed in \cite{chung2017constraint},
which rely on minimal dimensional
coarse spaces as discussed above.
First, we construct local auxiliary basis functions following the minimal dimensional coarse spaces
as discussed above.
For each coarse-block $K\in \mathcal{T}^H$, we solve the eigenvalue problem with Rayleigh quotient
$\mathcal{Q}_{ms}(v):=\frac{\int_{K}  \kappa  |\nabla v|^2}{\int_{K} \widehat{\kappa}| v  |^2}$,
where $\widehat{\kappa}=\kappa\sum_{j}|\nabla\chi_{j}|^{2}$. We assume
$\lambda^{K}_{1}\leq\lambda^{K}_{2}\leq\dots$ and define the local
auxiliary  spaces by
$$
V_{aux}(K)=\text{span}\{\phi_{j}^{K}|1\leq j\leq L_K\} \mbox{ and } V_{aux}=\oplus_{K}V_{aux}(K).$$
Next, define a projection operator $\pi_K$ as the orthogonal projection on $V_{aux}$ with respect to the inner product $\int_K \widehat{\kappa} uv$
and $\pi_D= \oplus_{K}\pi_K$.

Let $K^+$ be obtained by adding $l$ layers of coarse elements to the coarse-block $K$.
The coarse-grid multiscale basis  $\psi_{j,ms}^{K}\in V(K^{+})=
P^1_0(K^{+})$ solve
\[
\int_{K^+}\!\!\kappa\nabla \psi_{j,ms}^{K} \nabla v+\int_{K^+}\!\!\!\!
\widehat{\kappa}\pi_D(\psi_{j,ms}^{K})\pi_D(v)=
\int_{K^+} \widehat{\kappa}\phi_{j}^{K} \pi_D (v),\;
\forall v\in V(K^{+}).
\]
The coarse-grid multiscale  space is defined as $V_{ms}=\text{span}\{\psi_{j,ms}^{(i)}\}.$

Before discussing the method  using this coarse-grid space,
we  introduce some operators.
We consider the (coarse solution) operator $A_{0,ms}^{-1}:L^{2}(\Omega)\mapsto V_{ms}$ by
\[
\int_D \nabla A_{0,ms}^{-1}(u) \nabla v=\int_{\Omega}uv\;\;\mbox{ for all } v\in V_{ms}
\]
and the (local  solutions) operators $A_{i,ms}^{-1}:L^{2}(\Omega)\mapsto V(\omega_{i}^{+})$ defined by,
\[
\int_{\omega_i}\kappa \nabla A_{i,ms}^{-1}(u_{i}) \nabla v+
\int_{\omega_i} \widehat{\kappa}\pi(A_{i}^{-1}(u_{i})) \pi(v)=\int_{\omega_{i}}\chi_{i}uv\;\;\mbox{ for all } v\in P^1(\omega_{i}^{+}),
\]
where $\omega_{i}^{+}$ is obtained by enlarging
$\omega_{i}$ by $k$ coarse-grid layers. Next, we can define the preconditioner\footnote{Here we avoid restriction and extension operators for simplicity} $M$ by
\[
M^{-1}=(I-A_{0,ms}^{-1}A)\Big(\sum_{i}A_{i,ms}^{-1}\Big)(I-AA_{0,ms}^{-1})+A_{0,ms}^{-1}.
\]
Note that this is a hybrid preconditioner as defined in \cite{tw}.
Using some estimates in \cite{chung2017constraint},
we can show the bound of the form,
\begin{equation}\label{b2}
\text{cond}(M^{-1}A)\leq\cfrac{1+C(1+\Lambda^{-1})^{\frac{1}{2}}E^{\frac{1}{2}}\max\{\tilde{\kappa}^{\frac{1}{2}}\}}{1-C(1+\Lambda^{-1})^{\frac{1}{2}}E^{\frac{1}{2}}\max\{\tilde{\kappa}^{\frac{1}{2}}\}}
\end{equation}
where $E = 3(1+\Lambda^{-1}) \Big (1+(2(1+\Lambda^{-\frac{1}{2}}))^{-1}\Big)^{1-k}$, $C$ is a constant depend on the fine and coarse grid only and $\Lambda = \min_{K}{\lambda^{K}_{L_K+1}}$. See \cite{chung2017constraint}  for the required estimates of the coarse space. The analysis of the local solvers of the hybrid method above will be presented elsewhere due to the page limitation.
Here, we metion that the analysis do not use a stable decomposition so, in principle, a new family of robust method can be obtained. Moreover, we see
that the condition number is close to $1$ if sufficient number
of basis functions are selected (i.e., $\Lambda$ is not close to zero)\footnote{Having robust condition number close to 1 is important, specially in applications where the elliptic equation needs to be solved many times.}.
The overlap size usually involves several coarse-grid block sizes
and thus, the method is effective when the coarse-grid sizes are small.
We comment that taking the generous overlap $\delta=kH/2$ in \eqref{b1}, we get the bound $C(1+4/k^2)$  with $C$ independent of the contrast. The estimate
(\ref{b2}), on the other hand, gives a bound close to 1 if the oversampling
is sufficiently large (e.g., the number of coarse-grid layers is related
to $\log(\eta)$), which is due to the localization of global fields orthogonal to the coarse space.

\begin{figure}[htb]
\centering
\includegraphics[scale=0.4]{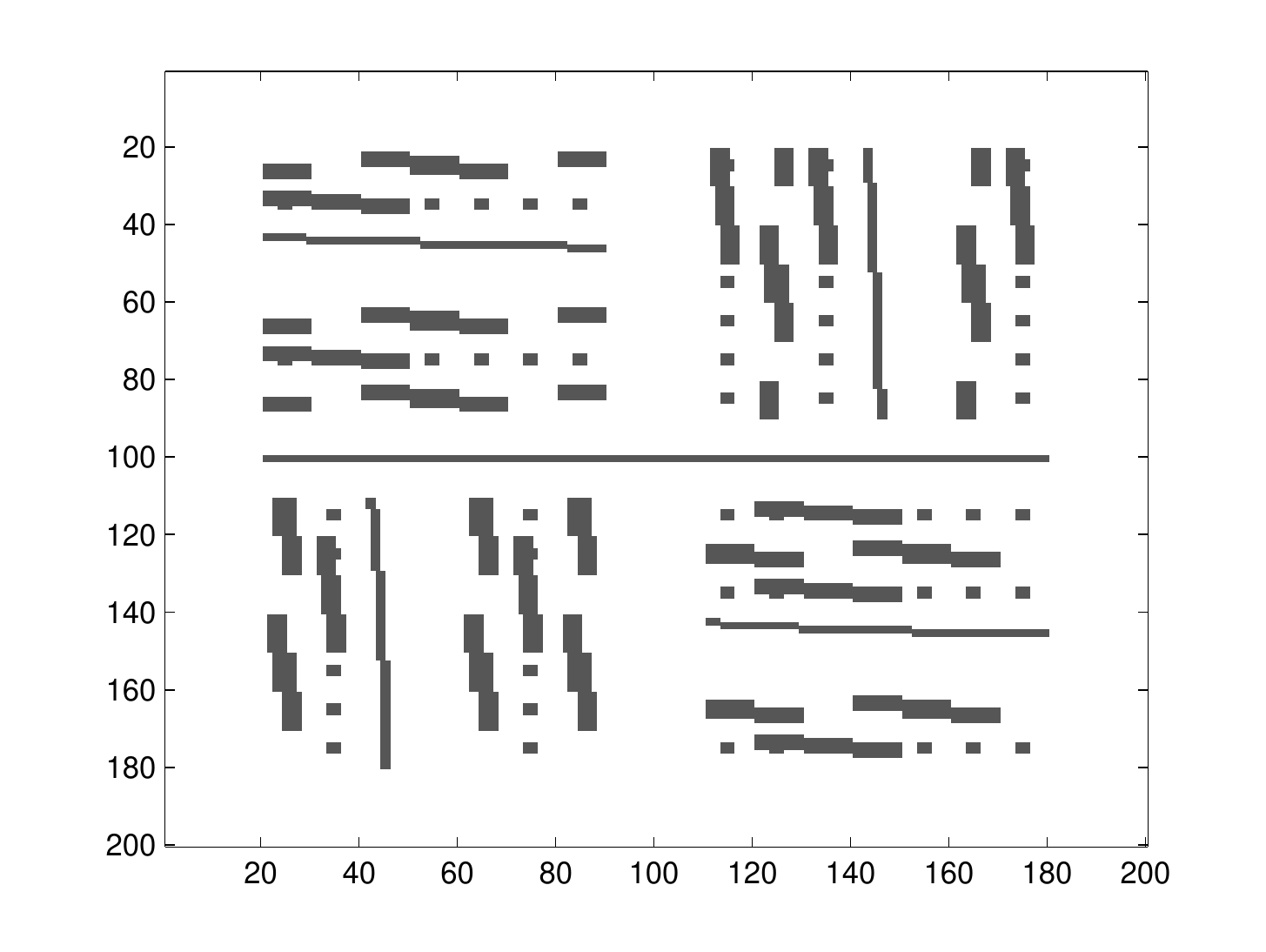}
\caption{The permeability $\kappa$ used in the numerical experiments. The grey regions indicate high-permeability region of order $\eta$ while the white regions indicates a low (order 1) permeability.}
\label{fig:medium}
\end{figure}

Next, we present a numerical result and consider a problem with permeability $\kappa$ shown in Fig.~\ref{fig:medium}. The fine-grid mesh size $h$ and the coarse-grid mesh size are considered as $h=1/200$ and $H=1/20$. In Table \ref{tab:result_CEM}, we present the number of iterations for using varying number of oversampling layers $k$ and value of
the contrast $\eta$.

\begin{figure}[htb]
\centering
\includegraphics[scale=0.3]{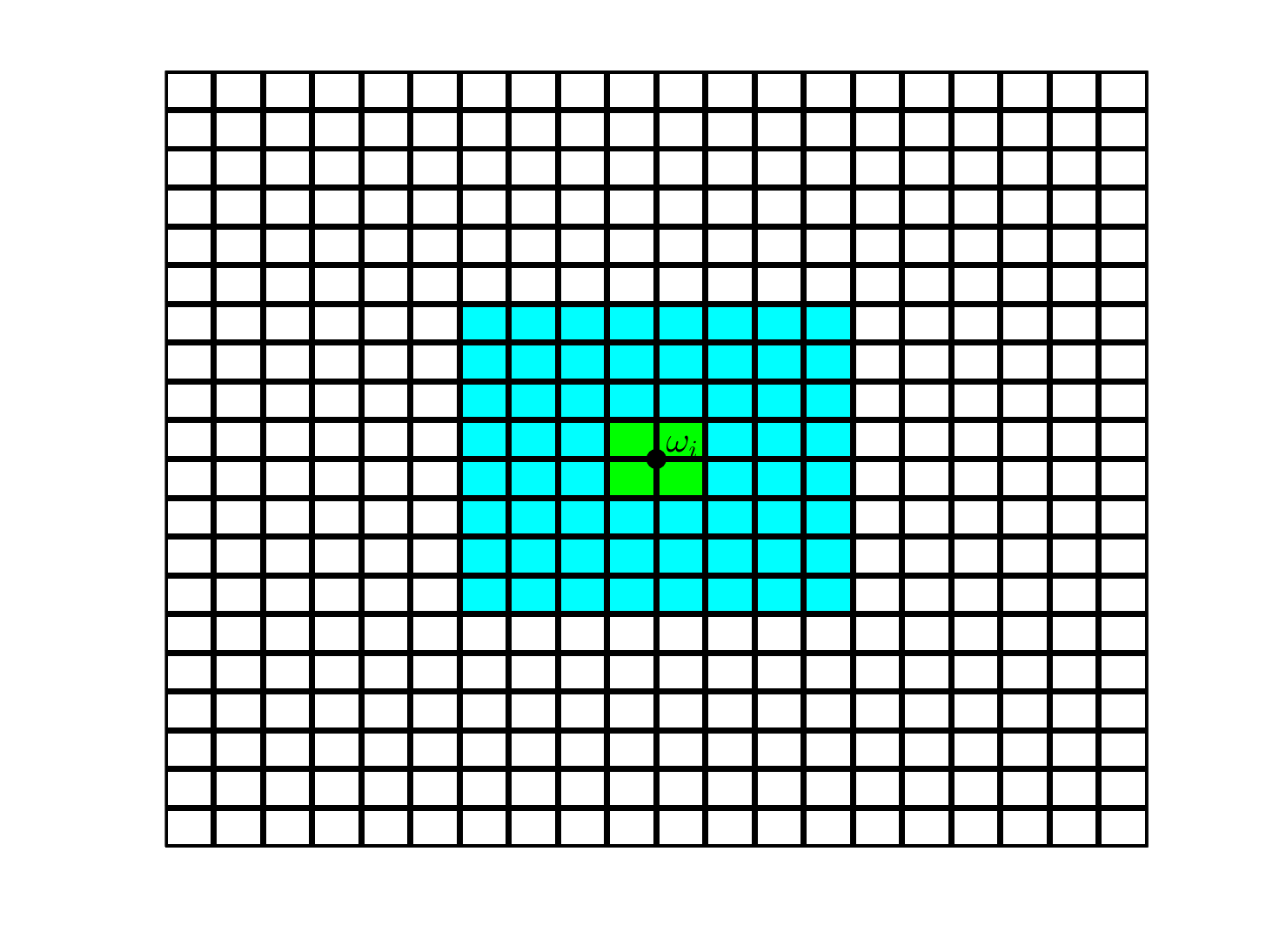}
\caption{The coarse mesh used in the numerical experiments. We highlight a coarse neighborhood and the results of adding 3 coarse-block layers to this neighborhood.}
\label{fig:coarse_mesh}
\end{figure}

\begin{table}[htb]
\centering
\begin{tabular}{|c|r|r|}
\hline
Number basis per $\omega$ & k & \# iter\tabularnewline
\hline
3  & 3 & 3\tabularnewline
3  & 4 & 2\tabularnewline
3  & 5 & 2\tabularnewline
3  & 6 & 1\tabularnewline
\hline\hline\hline
\end{tabular} \hspace{3pt}
\begin{tabular}{|c|c|r|}
\hline
Number basis per $\omega$ & $\eta$ & \# iter\tabularnewline
\hline
3 & 1e+3 & 3\tabularnewline
3 & 1e+4 & 3\tabularnewline
3 & 1e+5 & 3\tabularnewline
\hline\hline\hline
\end{tabular}

\caption{Number of iterations until convergence  for  the PCG with $H=1/20$, $h=1/200$ and $\text{tol}=1e-8$. Left: different number of oversampleing layers $k$ with $\eta=1e+4$. Right: different values of
the contrast $\eta$ with $k=3$.}
\label{tab:result_CEM}
\end{table}

We would like to emphasize that the proposed
method has advantages if the coarse mesh size is not very coarse.
In this case, the oversampled coarse regions are still sufficiently
small and the coarse-grid solves can be relatively expensive.
Consequently, one wants to minimize the number of coarse-grid solves
in addition to local solves.
In general, the proposed approach can be used in a multi-level setup,
in particular, at the finest levels, while at the coarsest level,
we can use original spectral basis functions proposed in \cite{ge09_1}. This is object of future research.

\vspace{-0.5cm}
\section{Conclusions}
\vspace{-0.5cm}
In this paper, we give an overview of domain decomposition preconditioners for multiscale high-contrast problems.
We emphasize the use of minimal dimensional coarse spaces in order to construct
optimal preconditioners with the condition number independent of physical scales (contrast and
spatial scales). We discuss various approaches in this direction. Furthermore, using these spaces
and oversampling ideas, we design a new preconditioner with significant reduction in the number of iterations until convergence if oversampling regions are large enough (several coarse-grid blocks).
We note that when using only minimal dimensional coarse spaces in additive Schwarz preconditioner with standard overlap size, we obtain around $19$ iterations. in the new method, our main goal is to reduce even further the number of iteration due to large coarse problem sizes. We obtained around 3 iteration until convergence for the new approach.
A main point of the new methodology is that after removing the channels we are able to localize the remaining multiscale information via oversampling.
Other interesting aspect of the new approach is that the bound can be obtained by estimating directly operator norms and do not require a stable decomposition.

\bibliographystyle{plain}
\bibliography{references,references1}

\end{document}